\magnification=\magstep1
\documentstyle{amsppt}
\PSAMSFonts
\hsize=6.5truein
\vsize=8.9truein
\NoRunningHeads
\refstyle{A}

\document

\topmatter
\baselineskip=26truept

\title 
Green's Theorem with No Differentiability 
\endtitle

\author 
Isidore Fleischer
\endauthor

\address
Centre de recherches  mathematiques, Universite de
Montreal, C.P. 6128, succ. Centre-ville, Montreal, QC H3C 3J7,Canada.
\endaddress

\abstract 

The result is established for a Jordan measurable region with rectifiable
boundary.  The integrand F for the new plane integral to be used is a
function of axis-parallel rectangles, finitely additive on non-overlapping
ones, hence unambiguously defined and additive on "figures" (i.e. finite
unions of axis-parallel rectangles).  Define its integral over Jordan
measurable S as the limit of its value on the figures, which contain a
subfigure of S and are contained in a figure containing S, as the
former/complements of the latter expand directedly to fill out S/the
complement of S.  The integral over every Jordan measurable region exists
when additive F is "absolutely continuous" in the sense of converging to
zero as the area enclosed by its argument does, or with F the
circumferential line integral  $\oint P\, dx + Q\, dy$  for $P$, $Q$
continuous at the rectifiable boundary of  S and
integrable along axis - parallel line segments.  Thus the equality of this
area integral with the line integral around the boundary, to be proved,
follows for the various integrals of divergence presented in:  Pfeffer,
W.F. The Riemann Approach to Integration, Cambridge Univ. Press, New York,
1993.    
\endabstract
\subjclass 26B20 26A39 28A75
\endsubjclass

\endtopmatter

In Advanced Calculus texts Green's Theorem is presented for continuous
vector fields with continuous first partial derivatives, defined in a
region containing a simple piecewise smooth curve enclosing an area of
not too complicated shape.  More careful treatments---e.g\. \cite{A,
\S10--14}---dispense with the continuity of the derivatives in favour
of their (bounded existence and) integrability over the interior;
recently this requirement has been successively weakened further to
integrability of the partials in the ``generalized Riemann'' sense
\cite{McL, \S7.12} and beyond to ``gage integrability'' \cite{Pf}
which even follows from the mere existence of the derivative.  By
modifying this last integral further, it proves possible to obtain the
theorem for a continuous vector field with no differentiability
assumption whatsoever.

\head
The ``integral'' of an additive rectangle function over a Jordan
measurable set
\endhead

For plane Jordan content, see \cite{A, \S10-4}, \cite{K}: Inner Jordan
content of a bounded set $S$ is the sup of areas of
finite unions of axis-parallel rectangles 
in the interior of
$S$; outer Jordan content, the inf of areas of such which meet closure of 
$S$. Their equality, ``Jordan measurability'', comes to the boundary
of $S$ having Jordan content (equivalently, by compactness, Lebesgue
measure) zero.

Let $F$ be a function of axis-parallel rectangles, finitely additive
on non-overlapping ones, hence unambiguously defined and additive on
``figures'' (i.e. finite unions of such).  We propose to define its 
integral over  $S$
as the common limit of its values, on the partition's 
(into axis-parallel rectangles of a rectangle containing $S$) 
smallest subfigure which contains the closure 
of $S$ or largest contained in its interior, as the partition is refined. 
The integral over every Jordan measurable
set $S$ exists when additive $F$ is ``absolutely continuous'' in the sense
of converging to zero as the area enclosed by its argument does---indeed,
the difference between a containing and a
contained figure of $S$ eventually has small area and so the value of
$F$ on the difference would be small, hence the values of
$F$ on these figures eventually close.  This obtains for the more 
usual kinds of area
integral of a bounded integrand. However, independent of this
absolute continuity,  the circumferential
line integral $\oint P\, dx + Q\, dy$, for $P$, $Q$ integrable along 
axis-parallel line segments and continuous at
the boundary (if rectifiable), construed as a 
rectangle function, will now be shown integrable in this sense 
to the usual line integral around the boundary as value. 

\head
Green's Theorem for Jordan regions with rectifiable boundaries
\endhead

Cover the boundary 
$\Gamma$ with non-overlapping axis-parallel rectangles of sufficiently 
small side length (finitely many by compactness) so that
each encloses only a connected piece of $\Gamma$ (which is locally
connected, since it is the continuous image of a compact such) and 
enters and leaves the rectangle through the shorter opposite sides 
(if part of $\Gamma$ is a horizontal or vertical interval, count this as a 
degenerate rectangle). Thus the sum of their perimeters is at most four
times the arc length of $\Gamma$. Construe the boundary rectangle sides in
the bounded  component of $S$ as the boundary of a figure in the interior of
$S$,  those in the unbounded component as the boundary of a cover.  

It is classical that the line integral of a continuous vector-valued
function along a rectifiable curve $\Gamma$ exists \cite{A, \S10-10}.
Since a constant vector line integrates around a closed curve to zero,
every  circumferential
line integral is bounded by the oscillation of the vector integrand
times the arc length of the circumference.
Adding up the line integrals (of a vector function continuous at $\Gamma$
and integrable over axis-parallel lines) around the circumferences (of
pieces of
$\Gamma$ and adjacent  rectangle sides) 
shows the line integral along $\Gamma$ 
close to that along each of the boundaries (if sufficiently close to
$\Gamma$) --- and of course along 
boundaries derived from any smaller rectangle cover of $\Gamma$, 
which thus bound larger interior figures and smaller covering figures of 
$S$ --- hence equal to the above defined surface integral for the 
$\lq\lq$flux": 
i.e. the circumferential line integral of the vector field construed as a
rectangle function.

\head
Supplement
\endhead

The same argument works in higher dimensions; 
I carry it through for dimension three   
(thus replacing $\lq\lq$Green" with $\lq\lq$Gauss").
Let $\Gamma$ be a homeomorph of the $2$-sphere, which bounds a Jordan
measurable volume $S$ in $3$-space and is ``rectifiable'': i.e.\ of  finite
total surface area and a.e.\ (for surface area) $C'$ --- then every
continuous vector-valued  function's normal component is integrable over
$\Gamma$. Cover $\Gamma$ with non-overlapping axis-parallel paralellepipeds
of sufficiently small  side area (finitely many by compactness) so that
each  meets $\Gamma$ in a connected piece and whose boundary meets $\Gamma$
in opposite sides of area no greater than 
that of the non-meeting pair (if part of $\Gamma$ is a horizontal or
vertical rectangle, count this as  a degenerate paralellepiped).
Thus the sum of their perimeters is at most six times the area of 
$\Gamma$. Construe the boundary rectangles in the bounded 
component of $S$ as the boundary of a figure in the interior of $S$, 
those in the unbounded component as the boundary of a cover.  
Since a constant vector's normal component
integrates around a closed surface to zero, every 
circumferential 
surface integral is bounded by the oscillation of the vector integrand
times the area of the perimeter.
Adding up the surface integrals (of a vector field's normal component, 
continuous at $\Gamma$  and integrable over axis-parallel rectangles) 
around the circumferences (of pieces of $\Gamma$ and adjacent 
rectangles) 
shows the surface integral of the normal component
along $\Gamma$ 
close to that along each of the boundaries (if sufficiently close to
$\Gamma$) --- and of course along 
boundaries derived from any smaller paralellepiped cover of $\Gamma$, 
which thus bounds larger interior figures and smaller covering figures of 
$S$ 
--- hence equal to the volume integral for the $\lq\lq$flux": 
i.e. the circumferential surface integral of the vector field construed as a
function of  axis-parallel paralellepipeds.

This has applications to the integrals discussed in [PF], which all 
feature additive rectangle functions under various supplementary 
hypotheses
of differentiability. His ``solids'' are Jordan measurable and so his
``divergence theorems'' would 
receive a better specified formulation from the above.

\Refs
\widestnumber\key{McL}

\vskip.10in
\parindent=0pt
\hangindent=15pt
\hangafter=1
[A] \quad T.M. Apostol, {\it Mathematical Analysis}, Addison-Wesley, 
Reading, MA, 1957.

\vskip.10in
\parindent=0pt
\hangindent=15pt
\hangafter=1
[G] \quad H. Grunsk L.Y., {\it The General Stokes' theorem}, Pitman, 1983.

\vskip.10in
\parindent=0pt
\hangindent=15pt
\hangafter=1
[K] \quad M.I. Knopp, {\it Theory of Area}, Markham, Chicago, 1969.

\vskip.10in
\parindent=0pt
\hangindent=15pt
\hangafter=1
[M] \quad F. Morgan, {\it Geometric Measure Theory}, Academic Press, 1988.

\vskip.10in
\parindent=0pt
\hangindent=15pt
\hangafter=1
[McL] \quad R.M. McLeod, {\it The Generalized Riemann Integral}, Carus, 
Math. Monograph, vol 20, 1980.

\vskip.10in
\parindent=0pt
\hangindent=15pt
\hangafter=1
[Pf] \quad W.F. Pfeffer, {\it The Riemann Approach to Integration}, 
Cambridge, Univ. Press, New York, 1993.

\vskip.10in
\parindent=0pt
\hangindent=15pt
\hangafter=1
[S] \quad S. Saks, {\it Theory of the Integral}, Hafner, New York, 1937.

\vskip.10in
\parindent=0pt
\hangindent=15pt
\hangafter=1
[SM] \quad K.T. Smith, {\it A Primer of Modern Analysis}, 
Springer Verlag, New York, 1983.

\endRefs
\enddocument